\documentclass[11pt]{article}
\usepackage[left=1in,right=1in,top=1in,bottom=1in]{geometry}
\usepackage{times}
\usepackage{expl3}
\usepackage{cite}
\usepackage[table]{xcolor}
\usepackage{multirow}
\usepackage{stackengine} 
\usepackage{hhline}
\usepackage{lipsum}
\usepackage{titlesec}
\usepackage{wrapfig}
\usepackage{enumerate}
\usepackage{epsfig}
\usepackage{enumitem}
\usepackage{bbm}
\usepackage{calc}
\usepackage{graphicx}
\usepackage{amsmath}
\usepackage[title]{appendix}
\usepackage{amssymb}
\usepackage{epstopdf}
\usepackage{boldline}
\usepackage{calligra}
\usepackage{bm}
\usepackage{url}
\usepackage{blindtext}
\usepackage{accents}

\newcommand{\define}{\stackrel{\mbox{\tiny def}}{=}}

\newtheorem{definition}{Definition}
\newtheorem{theorem}{Theorem}

\newtheorem{lemma}{Lemma}

\usepackage{mathtools}
\usepackage{epstopdf}
\usepackage{balance}
\usepackage{thmtools}
\usepackage{thm-restate}
\usepackage{hyperref}
\usepackage{cleveref}
\usepackage[mathscr]{euscript}

\usepackage[ruled,vlined]{algorithm2e}
\include{pythonlisting}
\newcommand{\ostar}{\mathbin{\mathpalette\make@circled\star}}

\makeatletter
\newcommand{\removelatexerror}{\let\@latex@error\@gobble}
\makeatother
\makeatletter
\newcommand*{\rom}[1]{\expandafter\@slowromancap\romannumeral #1@}
\makeatother

\ExplSyntaxOn
\newcommand\latinabbrev[1]{
  \peek_meaning:NTF . {
    #1\@}%
  { \peek_catcode:NTF a {
      #1.\@ }%
    {#1.\@}}}
\ExplSyntaxOff

\titleclass{\subsubsubsection}{straight}[\subsubsection]

\begin{document}
\vspace{1cm}
\title{Algebraic Connectivity Characterization of Ensemble Random Hypergraphs}\vspace{1.8cm}
\author{Shih~Yu~Chang
\thanks{Shih Yu Chang is with the Department of Applied Data Science,
San Jose State University, San Jose, CA, U. S. A. (e-mail: {\tt
shihyu.chang@sjsu.edu}). 
           }}

\maketitle

\begin{abstract}
Random hypergraph is a broad concept used to describe probability distributions over hypergraphs, which are mathematical structures with applications in various fields, e.g., complex systems in physics, computer science, social sciences, and network science. Ensemble methods, on the other hand, are crucial both in physics and machine learning. In physics, ensemble theory helps bridge the gap between the microscopic and macroscopic worlds, providing a statistical framework for understanding systems with a vast number of particles. In machine learning, ensemble methods are valuable because they improve predictive accuracy, reduce overfitting, lower prediction variance, mitigate bias, and capture complex relationships in data. However, there is limited research on applying ensemble methods to a set of random hypergraphs. This work aims to study the connectivity behavior of an ensemble of random hypergraphs. Specifically, it focuses on quantifying the random behavior of the algebraic connectivity of these ensembles through tail bounds. We utilize Laplacian tensors to represent these ensemble random hypergraphs and establish mathematical theorems, such as Courant-Fischer and Lieb-Seiringer theorems for tensors, to derive tail bounds for the algebraic connectivity. We derive three different tail bounds, i.e., Chernoff, Bennett, and Bernstein bounds, for the algebraic connectivity of ensemble hypergraphs with respect to different random hypergraphs assumptions.
\end{abstract}


\section{Introduction}\label{sec: Introduction}

Mathematics encompasses graph theory, a field dedicated to the examination of graphs—mathematical constructs employed for representing pairwise relationships among elements. In this mathematical context, a graph consists of vertices (referred to as nodes or points) that are interlinked by edges (also known as links or lines). To extend beyond conventional graphs, a natural progression involves permitting edges that can establish connections among an unrestricted number of nodes. A hypergraph is a type of graph where hyperedges, or generalized edges, have the capability to connect not just two vertices or nodes but rather a subset of vertices. In a hypergraph, these edges, often referred to as hyperedges, are defined as arbitrary nonempty sets of vertices.The utilization of hypergraphs is widespread in numerous fields in science and engineering~\cite{bunke2005theoretical,bretto2013applications}. 

The study of hypergraph connectivity is essential in various fields, including mathematics, computer science, and network analysis, due to its wide-ranging applications and significance~\cite{bapat2010graphs}. For instance,  connectivity problems help ensure reliable data transmission, fault tolerance, and efficient routing~\cite{gao2020v2vr}. In graph theory, the connectivity of a connected graph can be assessed through measures like edge connectivity and vertex connectivity. However, computing these measures can be notably challenging~\cite{chuzhoy2009k}. Fiedler introduced the concept of algebraic connectivity, establishing a relationship between algebraic connectivity and vertex connectivity~\cite{fiedler1973algebraic}. This connection allows for the estimation of vertex connectivity at a relatively lower computational cost. The notion of algebraic connectivity, as proposed by Fiedler, has also been extended to encompass directed graphs~\cite{wu2005algebraic} and is widely deliberated in the context of hypergraphs~\cite{asadi2016generalized}. In their work, Hu and Qi~\cite{hu2012algebraic} introduced a representation of a $2m$-uniform hypergraph using tensors. They employed $H$-eigenvalues and $Z$-eigenvalues as tools for defining the algebraic connectivity~\cite{li2013z}. In a related study, Cooper and Dutle delved into the exploration of spectral hypergraph theory by utilizing adjacency tensors~\cite{cooper2012spectra}.

The term \emph{random hypergraph} is a broad descriptor used to denote probability distributions across hypergraphs. Random hypergraphs are a fundamental concept in graph theory and have applications and significance in various fields, including physics, computer science, social sciences, and network analysis~\cite{ghoshal2009random}. Following examples show how random hypergraphs help us to solve problems in wide spectrum of areas. In information technology, random hypergraphs are used to model the internet's topology and computer networks. Understanding network properties can lead to better network design, routing algorithms, and resource allocation~\cite{alon1999independent}. In physics, random hypergraph structures exhibit phase transitions, where certain properties of the hypergraph change suddenly as parameters (such as the edge probability) vary. Understanding these phase transitions provides insights into the emergence of network structures in condense matters~\cite{ferraz2021phase}. In computer science random hypergraphs are used in the design and analysis of randomized algorithms. These algorithms use randomness to solve problems efficiently or approximate solutions to problems that are otherwise computationally challenging~\cite{fohlin2008randomized}. In sociology and social science, random hypergraph models are employed to study social networks, diffusion of information, and social influence. They help explain patterns of social connections and behavior~\cite{chodrow2019configuration}.

On the other hand, the concept of \emph{ensemble} is important in physics and learning theory. In physics, particularly in statistical mechanics, ensemble theory helps us understand and describe the behavior of large collections of particles or systems, such as gases, solids, and liquids. It is essential because it provides a statistical framework for dealing with systems that involve a vast number of particles, making it infeasible to track the behavior of each individual particle. The main reason to apply ensemble concepts to interpret physical phenomena is that ensemble theory helps scientists to bridge the gap between the microscopic and macroscopic worlds, making it possible to describe and predict the behavior of large collections of particles or systems. It is a foundational concept that underpins our understanding of thermodynamics, statistical mechanics, and many other areas of physics~\cite{van2015ensemble}. In learning theory, ensemble methods are a valuable approach based on the following benefits. First, ensemble methods often lead to better predictive performance compared to individual models. By combining multiple models, they can reduce errors and increase overall accuracy. Second, ensembles of learning models are less prone to overfitting, which occurs when a model learns to perform well on the training data but fails to generalize to new, unseen data. Ensembles can generalize better because they capture diverse patterns and reduce the risk of modeling noise. Third, ensemble methods can reduce the variance of predictions. This is particularly beneficial when dealing with noisy or uncertain data. Fourth, ensembles can mitigate bias present in individual models. By combining models with different biases, ensembles can provide a more balanced prediction. Finally, ensembles can capture complex relationships in the data by combining the outputs of multiple models with different perspectives or assumptions. This makes an ensemble approach to expand the space of representable functions for solutions~\cite{dietterich2000ensemble}. However, there are not many works about considering ensemble behavior for a set of random hypergraphs.

The purpose of this work is to study the connectivity behavior for ensembles of random hypergraphs. Given a set of random $2m$-uniform hypergraphs, represented by $\{\mathscr{G}_i\}_{i=1}^N$, we wish to quantify the random behavior of the algebraic connectivity for ensemble of these random hypergraphs via the tail bound of the algebraic connectivity which can be expressed by:
\begin{eqnarray}\label{eq:prob of connect formulation}
\mathrm{Pr}\left(\alpha(\overline{\mathscr{G}}) \geq \theta \right),
\end{eqnarray}
where $\alpha(\overline{\mathscr{G}})$ is the algebraic connectivity of the ensemble random hypergraph, represented by $\overline{\mathscr{G}}$ and $\theta$ is a specified positive number about the degree of connectivity. The Laplacian tensor of the ensemble random hypergraph $\overline{\mathscr{G}}$, denoted by $\mathcal{L}_{\overline{\mathscr{G}}}$, is obtained by the following relation:
\begin{eqnarray}\label{eq:def:ensemble Laplacian}
\mathcal{L}_{\overline{\mathscr{G}}} = \sum\limits_{i=1}^N \mathcal{L}_{\mathscr{G}_i},
\end{eqnarray}
where $\mathcal{L}_{\mathscr{G}_i}$ are Laplacia tensors for uniform random hypergraphs $\{\mathscr{G}_i\}$. The detailed notions about  $2m$-uniform hypergraphs, algebraic connectivity, Laplacian tensors are discussed in Section~\ref{sec: Algebraic Connectivity of A Hypergraph}. According to~\cite{gu2023even}, the algebraic connectivity of the ensemble random hypergraph is the $k$-th largest eigenvalue of the Laplcian tensor $\mathcal{L}_{\overline{\mathscr{G}}}$, where $k$ is some positive integer related to hypergraph parameters of $\overline{\mathscr{G}}$. Therefore, in order to bound the probability provided by Eq.~\eqref{eq:prob of connect formulation}, it is equivalent to determine the tail bound for the $k$-th largest eigenvalue for the summation of random tensors. We establish Courant-Fischer theorem and Lieb-Seiringer joint concavity theorem for tensors first, and apply these theorems to obtain Chernoff, Bennett, and Bernstein tail bounds for the algebraic connectivity $\overline{\mathscr{G}}$ under different random hypergraphs $\mathscr{G}_i$ assumptions. 

The rest of this paper is organized as follows. Basic tensor facts and functions of tensors are introduced in Section~\ref{sec: Tensor Review}. In Section~\ref{sec: Algebraic Connectivity of A Hypergraph}, we will discuss the random hypergraph and its algebraic connectivity. In Section~\ref{sec: Tail Bounds for $k$-th Largest Eigenvalue}, we will prepare the tail bounds for any $k$-th largest eigenvalue arising from the summation of independent Hermitian random tensors. Finally, in Section~\ref{sec: Algebraic Connectivity of Ensemble Hypergraphs}, we will consider different tail bounds, i.e., Chernoff, Bennett, and Bernstein bounds, for the algebraic connectivity of ensemble hypergraphs with respect to different random hypergraphs assumptions. 

\textbf{nomenclasure}

Throughout this work, tensors are represented by calligraphic letters (e.g., $\mathcal{D}$, $\mathcal{E}$, $\mathcal{F}$, $\ldots$), respectively. Tensors are multiarrays of values which are higher-dimensional generalizations from vectors and matrices. Given a positive integer $N$, let $[N] \define \{1, 2, \cdots ,N\}$. An \emph{order-$M$ tensor} (or \emph{$M$-th order tensor}) denoted by $\mathcal{X} \define (a_{i_1, i_2, \cdots, i_M})$, where $1 \leq i_j = 1, 2, \ldots, I_j$ for $j \in [M]$, is a multidimensional array containing $I_1 \times I_2 \times \cdots \times I_{M}$ entries. Let $\mathbb{C}^{I_1 \times \cdots \times I_M}$ and $\mathbb{R}^{I_1 \times \cdots \times I_M}$ be the sets of the order-$M$ $I_1 \times \cdots \times I_M$ tensors over the complex field $\mathbb{C}$ and the real field $\mathbb{R}$, respectively. For example, $\mathcal{X} \in \mathbb{C}^{I_1 \times \cdots \times I_M}$ is an order-$N$ multiarray, where the first, second, ..., and $M$-th dimensions have $I_1$, $I_2$, $\ldots$, and $I_M$ entries, respectively. Thus, each entry of $\mathcal{X}$ can be represented by $a_{i_1, \cdots, i_M}$. 

\section{Tensor Review}\label{sec: Tensor Review}

In this section, we'll offer a concise overview of tensors and probability, essential prerequisites for our forthcoming theory elaboration. This section is based on our previous work in~\cite{chang2022convenient}.

\subsection{Fundamental of Tensor}\label{subsec: Fundamental of Tensors}

\subsubsection{Tensor Notations}\label{sec:Tensor Notations}

We can divide the dimensions of a tensor into two distinct sets, denoted as $M$ and $N$ dimensions, respectively. Thus, for two tensors of order $M+N$: $\mathcal{X} \define (a_{i_1, \cdots, i_M, j_1, \cdots,j_N}) \in \mathbb{C}^{I_1 \times \cdots \times I_M\times
J_1 \times \cdots \times J_N}$ and $\mathcal{Y} \define (b_{i_1, \cdots, i_M, j_1, \cdots,j_N}) \in \mathbb{C}^{I_1 \times \cdots \times I_M\times
J_1 \times \cdots \times J_N}$, as described in~\cite{chang2022convenient}, the operation of \emph{tensor addition} $\mathcal{X} + \mathcal{Y}\in \mathbb{C}^{I_1 \times \cdots \times I_M\times
J_1 \times \cdots \times J_N}$ can be defined as follows:
\begin{eqnarray}\label{eq: tensor addition definition}
(\mathcal{X} + \mathcal{Y} )_{i_1, \cdots, i_M, j_1 \times \cdots \times j_N} &\define&
 a_{i_1, \cdots, i_M, j_1 \times \cdots \times j_N} \nonumber \\
& &+ b_{i_1, \cdots, i_M, j_1 \times \cdots \times j_N}. 
\end{eqnarray}
Conversely, when considering tensors $\mathcal{X} \define (a_{i_1, \cdots, i_M, j_1, \cdots, j_N}) \in \mathbb{C}^{I_1 \times \cdots \times I_M\times
J_1 \times \cdots \times J_N}$ and \\ $\mathcal{Y} \define (b_{j_1, \cdots, j_N, k_1, \cdots, k_L}) \in \mathbb{C}^{J_1 \times \cdots \times J_N\times K_1 \times \cdots \times K_L}$, as elucidated in~\cite{chang2022convenient}, the operation known as the \emph{Einstein product} (referred to simply as the \emph{tensor product} in this context) denoted as $\mathcal{X} \star_{N} \mathcal{Y} \in \mathbb{C}^{I_1 \times \cdots \times I_M\times
K_1 \times \cdots \times K_L}$ can be defined as follows:
\begin{eqnarray}\label{eq: Einstein product definition}
\lefteqn{(\mathcal{X} \star_{N} \mathcal{Y} )_{i_1, \cdots, i_M,k_1 \times \cdots \times k_L} \define} \nonumber \\ &&\sum\limits_{j_1, \cdots, j_N} a_{i_1, \cdots, i_M, j_1, \cdots,j_N}b_{j_1, \cdots, j_N, k_1, \cdots,k_L}. 
\end{eqnarray}
Please note that we will frequently condense the tensor product $\mathcal{X} \star_{N} \mathcal{Y}$ to "$\mathcal{X} \hspace{0.05cm}\mathcal{Y}$" for the sake of brevity in the remainder of this document. This tensor product simplifies to standard matrix multiplication when $L = M = N = 1$. Other simplified scenarios include the tensor-vector product ($M > 1$, $N = 1$, and $L = 0$) and the tensor-matrix product ($M > 1$ and $N = L = 1$). Similar to matrix analysis, we introduce fundamental tensors and elementary tensor operations as follows.

\begin{definition}\label{def: zero tensor}
A tensor whose entries are all zero is called a \emph{zero tensor}, denoted by $\mathcal{O}$. 
\end{definition}
\begin{definition}\label{def: identity tensor}
An \emph{identity tensor} $\mathcal{I} \in  \mathbb{C}^{I_1 \times \cdots \times I_N\times
J_1 \times \cdots \times J_N}$ is defined by 
\begin{eqnarray}\label{eq: identity tensor definition}
(\mathcal{I})_{i_1 \times \cdots \times i_N\times
j_1 \times \cdots \times j_N} \define \prod_{k = 1}^{N} \delta_{i_k, j_k},
\end{eqnarray}
where $\delta_{i_k, j_k} \define 1$ if $i_k  = j_k$; otherwise $\delta_{i_k, j_k} \define 0$.
\end{definition}
To establish the concept of a \emph{Hermitian} tensor, we outline the operation known as the \emph{conjugate transpose} (or \emph{Hermitian adjoint}) of a tensor in the following manner.
\begin{definition}\label{def: tensor conjugate transpose}
Given a tensor $\mathcal{X} \define (a_{i_1, \cdots, i_M, j_1, \cdots,j_N}) \in \mathbb{C}^{I_1 \times \cdots \times I_M\times J_1 \times \cdots \times J_N}$, its conjugate transpose, denoted by
$\mathcal{X}^{H}$, is defined by
\begin{eqnarray}\label{eq:tensor conjugate transpose definition}
(\mathcal{X}^H)_{ j_1, \cdots,j_N,i_1, \cdots, i_M}  \define  
\overline{a_{i_1, \cdots, i_M,j_1, \cdots,j_N}},
\end{eqnarray}
where the overline notion indicates the complex conjugate of the complex number $a_{i_1, \cdots, i_M,j_1, \cdots,j_N}$. If a tensor $\mathcal{X}$ satisfies $ \mathcal{X}^H = \mathcal{X}$, then $\mathcal{X}$ is a \emph{Hermitian tensor}. 
\end{definition}
\begin{definition}\label{def: unitary tensor}
Given a tensor $\mathcal{U} \define (u_{i_1, \cdots, i_M, i_1, \cdots,i_M}) \in \mathbb{C}^{I_1 \times \cdots \times I_M\times I_1 \times \cdots \times I_M}$, if
\begin{eqnarray}\label{eq:unitary tensor definition}
\mathcal{U}^H \star_M \mathcal{U} = \mathcal{U} \star_M \mathcal{U}^H = \mathcal{I} \in \mathbb{C}^{I_1 \times \cdots \times I_M\times I_1 \times \cdots \times I_M},
\end{eqnarray}
then $\mathcal{U}$ is a \emph{unitary tensor}. 
\end{definition}
In this work, the symbol $\mathcal{U}$ is resrved for a unitary tensor. 

\begin{definition}\label{def: inverse of a tensor}
Given a \emph{square tensor} $\mathcal{X} \define (a_{i_1, \cdots, i_M, j_1, \cdots,j_M}) \in \mathbb{C}^{I_1 \times \cdots \times I_M\times I_1 \times \cdots \times I_M}$, if there exists $\mathcal{X} \in \mathbb{C}^{I_1 \times \cdots \times I_M\times I_1 \times \cdots \times I_M}$ such that 
\begin{eqnarray}\label{eq:tensor invertible definition}
\mathcal{X} \star_M \mathcal{X} = \mathcal{X} \star_M \mathcal{X} = \mathcal{I},
\end{eqnarray}
then $\mathcal{X}$ is the \emph{inverse} of $\mathcal{X}$. We usually write $\mathcal{X} \define \mathcal{X}^{-1}$ thereby. 
\end{definition}

We also list other crucial tensor operations here. The \emph{trace} of a square tensor is equivalent to the summation of all diagonal entries such that 
\begin{eqnarray}\label{eq: tensor trace def}
\mathrm{Tr}(\mathcal{X}) \define \sum\limits_{1 \leq i_j \leq I_j,\hspace{0.05cm}j \in [M]} \mathcal{X}_{i_1, \cdots, i_M,i_1, \cdots, i_M}.
\end{eqnarray}
The \emph{inner product} of two tensors $\mathcal{X}$, $\mathcal{Y} \in \mathbb{C}^{I_1 \times \cdots \times I_M\times J_1 \times \cdots \times J_N}$ is given by 
\begin{eqnarray}\label{eq: tensor inner product def}
\langle \mathcal{X}, \mathcal{Y} \rangle \define \mathrm{Tr}\left(\mathcal{X}^H \star_M \mathcal{Y}\right).
\end{eqnarray}
According to Eq.~\eqref{eq: tensor inner product def}, the \emph{Frobenius norm} of a tensor $\mathcal{X}$ is defined by 
\begin{eqnarray}\label{eq:Frobenius norm}
\left\Vert \mathcal{X} \right\Vert \define \sqrt{\langle \mathcal{X}, \mathcal{X} \rangle}.
\end{eqnarray}

We use $\lambda_{\min}$ and $\lambda_{\max}$ to repsent the minimum and the maximum eigenvales of a Hermitain tensor. The notation $\succeq$ is used to indicate the semidefinite ordering of tensors. If we have $\mathcal{X} \succeq \mathcal{Y}$, this means that the difference tensor $\mathcal{X} - \mathcal{Y}$ is a positive semidefinite tensor.

\subsubsection{Tensor Functions}\label{sec:Tensor Functions}

For a tensor $\mathcal{A} \in \mathbb{C}^{I_1 \times \cdots \times I_M\times I_1 \times \cdots \times I_M}$, its eigenvalue $\lambda$ and corresponding eigen-tensor $\mathcal{X} \in \mathbb{C}^{I_1 \times \cdots \times I_M}$ is defined as
\begin{eqnarray}
\mathcal{A} \star_M \mathcal{X}&=& \lambda \mathcal{X}.
\end{eqnarray}
Given a function $g: \mathbb{R} \rightarrow \mathbb{R}$, the mapping result of a diagonal tensor by the function $g$ is to obtain another same size diagonal tensor with diagonal entry mapped by the function $g$. Then, the function $g$ can be extended to allow a Hermitian tensor $\mathcal{X} \in \mathbb{C}^{I_1 \times \cdots \times I_M\times I_1 \times \cdots \times I_M}$ as an input argument as
\begin{eqnarray}\label{eq:tensor func def}
g(\mathcal{X}) \define \mathcal{U} \star_M g(\Lambda) \star_M \mathcal{U}^H,~~\mbox{where $\mathcal{X} =  \mathcal{U}\star_M \Lambda \star_M \mathcal{U}^H $.}
\end{eqnarray}
The $\Lambda$ is the diagonal tensor with entries as eigenvalues. 

The \emph{spectral mapping theorem} posits that every eigenvalue of $g(\mathcal{X})$ corresponds to $g(\lambda)$ for some eigenvalue $\lambda$ of $\mathcal{X}$. Additionally, based on the semidefinite ordering of tensors, we obtain the following relationship:
\begin{eqnarray}\label{eq:tensor psd ordering}
f(x) \geq g(x),~~\mbox{for $x \in [a, b]$}~~~ \Rightarrow ~~~ f(\mathcal{X})  \succeq 
 g(\mathcal{X}),~~\mbox{for eigenvalues of $\mathcal{X} \in [a, b]$;}
\end{eqnarray}
where $[a, b]$ is a real interval. 

\begin{definition}\label{def: tensor exponential}
Given a square tensor $\mathcal{X} \in \mathbb{C}^{I_1 \times \cdots \times I_M\times I_1 \times \cdots \times I_M}$, the \emph{tensor exponential} of the tensor $\mathcal{X}$ is defined as 
\begin{eqnarray}\label{eq: tensor exp def}
e^{\mathcal{X}} \define \sum\limits_{k=0}^{\infty} \frac{\mathcal{X}^{k}}{k !}, 
\end{eqnarray}
where $\mathcal{X}^0$ is defined as the identity tensor $\mathcal{I} \in \mathbb{C}^{I_1 \times \cdots \times I_M\times I_1 \times \cdots \times I_M}$ and \\
$\mathcal{X}^{k} = \underbrace{\mathcal{X} \star_M \mathcal{X} \star_M \dots \star_M\mathcal{X} }_{\mbox{$k$ terms of $\mathcal{X}$}} $.

Given a tensor $\mathcal{Y}$, the tensor $\mathcal{X}$ is said to be a \emph{tensor logarithm} of $\mathcal{Y}$ if $e^{\mathcal{X}}  = \mathcal{Y}$
\end{definition}

Several important properties pertain to the tensor exponential. Firstly, as per the spectral mapping theorem, the exponential of a Hermitian tensor is always positive-definite. Secondly, the trace exponential function, denoted as $\mathcal{X} \rightarrow \mathrm{Tr} \exp(\mathcal{X})$, is convex. Thirdly, the trace exponential function exhibits a monotonic behavior with respect to the semidefinite ordering, as
\begin{eqnarray}\label{eq:trace exp monotone}
\mathcal{X} \succeq \mathcal{Y} ~~~ \Rightarrow ~~~ \mathrm{Tr} \exp(\mathcal{X}) \geq  \mathrm{Tr} \exp(\mathcal{Y}).
\end{eqnarray}

For the tensor logarithm, we have the following monotone relation 
\begin{eqnarray}\label{eq:log monotone}
\mathcal{X} \succeq \mathcal{Y} ~~~ \Rightarrow ~~~ \log(\mathcal{X}) \succeq  \log(\mathcal{Y}).
\end{eqnarray}
Moreover, the tensor logarithm is also concave, i.e., we have
\begin{eqnarray}\label{eq:log concave}
t \log (\mathcal{X}_1) + (1-t) \log (\mathcal{X}_2)  \preceq \log( t \mathcal{X}_1 +
(1-t) \mathcal{X}_2 ),
\end{eqnarray}
where $\mathcal{X}_1, \mathcal{X}_2$ are positive-definite tensors and $t \in [0,1]$. The concavity of tesnor logarithm can be dervied from Hansen-Pedersen Characterizations, see~\cite{MR3426646}

\subsection{Tensor Moments and Cumulants}\label{subsec: Tensor Moments and Cumulants}

Because the expectation of a random tensor can be treated as convex combination, expectation will preserve the semidefinite order as
\begin{eqnarray}
\mathcal{X} \succ \mathcal{Y} \mbox{~~almost surely} ~~~ \Rightarrow ~~~
\mathbb{E}(\mathcal{X}) \succ \mathbb{E}(\mathcal{Y}).
\end{eqnarray}
From operator Jensen's inequality~\cite{MR649196}, we also have 
\begin{eqnarray}
\mathbb{E}(\mathcal{X}^2) \succeq \left(\mathbb{E}(\mathcal{X})\right)^2.
\end{eqnarray}

Assume we have a random Hermitian tensor $\mathcal{X}$ with existing tensor moments of all orders, denoted as $\mathbb{E}(\mathcal{X}^n)$ for all $n$. In such a case, we can introduce the tensor moment-generating function, represented as $\mathbb{M}{\mathcal{X}}(t)$, and the tensor cumulant-generating function, represented as $\mathbb{K}_{\mathcal{X}}(t)$, for the tensor $\mathcal{X}$ as follows:
\begin{eqnarray}\label{eq:def mgf and cgf}
\mathbb{M}_{\mathcal{X}}(t) \define \mathbb{E} e^{t \mathcal{X}}, \mbox{~~and~~~}
\mathbb{K}_{\mathcal{X}}(t) \define \log \mathbb{E} e^{t \mathcal{X}},
\end{eqnarray}
where $t \in \mathbb{R}$. Both the tensor moment-generating function and the tensor cumulant-generating function can be expressedby power series expansions:
\begin{eqnarray}\label{eq:mgf and cgf expans}
\mathbb{M}_{\mathcal{X}}(t) = \mathcal{I} + \sum\limits_{n=1}^{\infty} \frac{t^n}{n !}\mathbb{E}(\mathcal{X}^n), \mbox{~~and~~~}
\mathbb{K}_{\mathcal{X}}(t) = \sum\limits_{n=1}^{\infty} \frac{t^n}{n !} \psi_n,
\end{eqnarray}
where $\psi_n$ is named as \emph{tensor cumulant}. The tensor cumulant $\psi_n$ can be expressed as a polyomial in terms of tensor moments up to the order $n$, for example, the first cumulant is the mean and the second cumulant is the varaince: 
\begin{eqnarray}\label{eq:cumulant mean and var}
\psi_1 = \mathbb{E}(\mathcal{X}), \mbox{~~and~~~}
\psi_2 = \mathbb{E}(\mathcal{X}^2) - (\mathbb{E}(\mathcal{X}))^2.
\end{eqnarray}

\section{Algebraic Connectivity of A Hypergraph}\label{sec: Algebraic Connectivity of A Hypergraph}

In this section, we first introduce the random hypergraph concept in Section~\ref{subsec: Random Hypergraph} and discuss algebraic connectivity of a hypergraph in Section~\ref{subsec: Algebraic Connectivity}.

\subsection{Random Hypergraph}\label{subsec: Random Hypergraph}

A hypergraph $\mathscr{G}=(\mathscr{V}, \mathscr{E})$ is composed by two parts: the set of vertices represented by $\mathscr{V}=\{V_1, V_2, \cdots, V_m\}$, and set of hyperedges represented by $\mathscr{E}=\{E_1, E_2, \cdots, E_n\}$. A hyperedge $E_i$ for $i \in [n]$ is a subset of $\mathscr{V}$. A hypergraph is named as $M$-uniform hypergraph if all its hyperedges contain the same number of vertices, i.e., $\left\vert E_i \right\vert =M$ for $i \in [n]$. A conventional graph with all edges determined by two vertices can be treated as a $2$-uniform hypergraph. In this work, we will consider only $2M$-uniform hypergraphs.  

According to the order among vertices in each $E_i$, we can have two types of hypergraph: undirected hypergraph and directd hypergraph. For an undirected hypergraph, there are no order relationships among vertices in each $E_i$; for a directd hypergraph, the vertices in $E_i$ are devided into two parts: the source part, denoted by $S_i$, and the destination part, denoted by $D_i$. $E_i = S_i \bigcup D_i$, and $S_i \bigcap D_i = \emptyset$. We have a weight function associated to each directed $\mathscr{G}=(\mathscr{V}, \mathscr{E})$ by assigning each $E_i$ a real value, denoted by $W_i$. A directed hypergraph $\mathscr{G}=(\mathscr{V}, \mathscr{E})$ with such weights assignments to each edge $E_i$ is named as a weighted hypergraph. For any two edges $E_i, E_j$ from a weighted $2M$-uniform hypergraph with $E_i = E_j$, $S_i = D_j$ and $D_i = S_j$, if $w_i = w_j$, such hypergraph is named as a \emph{weighted symmetric $2M$-uniform hypergraph}. If all weights $E_i$ of a weighted symmetric $2M$-uniform hypergraph are random variables, this hypergraph is called by \emph{random weighted symmetric $2M$-uniform hypergraph}. The proposed random weighted symmetric $2M$-uniform hypergraph is the extension of even order uniform hypergraphs (deterministic  hypergraphs) discussed in~\cite{gu2023even}. To the end of this work, we will focus on discussing random hypergraphs of this type: random weighted symmetric $2M$-uniform hypergraph.

For a random weighted symmetric $2M$-uniform hypergraph $\mathscr{G}=(\mathscr{V}, \mathscr{E})$ with $m$ vertices, we define its adjacency tensor, represented by $\mathcal{A}_{\mathscr{G}}=[a_{i_1,\cdots,i_m,j_1,\cdots,j_m}] \in \mathbb{R}^{m^{2M}}$, as follows:
\begin{eqnarray}\label{def:Adj tensor}
a_{i_1,\cdots,i_m,j_1,\cdots,j_m}&\define& w_{i_1,\cdots,i_m,j_1,\cdots,j_m},
\end{eqnarray}
where $w_{i_1,\cdots,i_m,j_1,\cdots,j_m}$ is the weight for edges between the source part with vertices $\{i_1,\cdots,i_m\}$ and the destination part with vertices $\{j_1,\cdots,j_m\}$. The degree tensor, represented by $\mathcal{D}_{\mathscr{G}}=[d_{i_1,\cdots,i_m,j_1,\cdots,j_m}] \in \mathbb{R}^{m^{2M}}$, is defined by
\begin{eqnarray}\label{def:Deg tensor}
d_{i_1,\cdots,i_m,j_1,\cdots,j_m}&\define& \sum\limits_{j_1,\cdots,j_M}^m a_{i_1,\cdots,i_m,j_1,\cdots,j_m}.
\end{eqnarray}
Then, the Laplacian tensor for this graph $\mathscr{G}$, represented by $\mathcal{L}_{\mathscr{G}}\in \mathbb{R}^{m^{2M}}$, is defined by 
\begin{eqnarray}\label{def:Lap tensor}
\mathcal{L}_{\mathscr{G}}&\define& \mathcal{D}_{\mathscr{G}} - \mathcal{A}_{\mathscr{G}}.
\end{eqnarray}

\subsection{Algebraic Connectivity of A Hypergraph}\label{subsec: Algebraic Connectivity}

Before defining the algebraic connectivity for a hypergraph, we say that edges $(E_1, E_2, \cdots, E_n)$ form a $M$-path if we have
\begin{eqnarray}
\left\vert E_i \bigcap E_{i+1} \right\vert = M,
\end{eqnarray}
where $i \in [n-1]$; and we also have
\begin{eqnarray}
E_i \bigcap E_{i+1} \bigcap E_{i+2} = \emptyset,
\end{eqnarray}
where $i \in [n-2]$.  We say that a given $\mathscr{G}$ is $M$-connected, i.e., for any $\{i_1, \cdots, i_M\}\subset \mathscr{V}$, and $\{j_1, \cdots, j_M\}\subset \mathscr{V}$, there exists an $M$-path from $\{i_1, \cdots, i_M\}$ to $\{j_1, \cdots, j_M\}$. 

From Eq. (2) in~\cite{gu2023even}, the algebraic connectivity for a hypergrah $\mathscr{G}$, denoted by $\alpha(\mathscr{G})$, is defined as follows:
\begin{eqnarray}\label{eq:alg def}
\alpha(\mathscr{G}) &\define& \min\limits_{\left\Vert \mathcal{X} \right\Vert,\mathcal{X}\perp \mathcal{E}}\mathcal{X}^{\mathrm{T}}\star_M \mathcal{L}_{\mathscr{G}} \star_M \mathcal{X},
\end{eqnarray}
where $\mathcal{E}$ is the subspace of $\mathbb{R}^{m^M}$ spanned by eigen-tensors
corresponding to the following two cases of zero eigenvalues:
\begin{enumerate}
\item All-one tensor with dimension $m^M \times 1$.
\item Within the adjacency tensor, duplicate indices result in a zero-filled section in the tensor $\mathcal{A}_{\mathscr{G}}$.
\end{enumerate}

According to Theorem 3.2 in~\cite{gu2023even} and the assumption that $m > M$, if a hypergraph $\mathscr{G}$ is a $M$-connected hypergraph, then we have
\begin{eqnarray}\label{eq: alg con connect cond}
\alpha(\mathscr{G}) \geq 0.
\end{eqnarray}
Moreover, the algebraic connectivity $\alpha(\mathscr{G})$ is the $\left(\frac{m!}{(m-M)!}-1\right)$-th largest eigenvalue of $\mathcal{L}_{\mathscr{G}}$ from Lemma 3.1 in~\cite{gu2023even}. 

\section{Tail Bounds for $k$-th Largest Eigenvalue}\label{sec: Tail Bounds for $k$-th Largest Eigenvalue}

In this section, we establish a universal upper limit for the tail probabilities of any $k$-th largest eigenvalue arising from the summation of independent Hermitian random tensors.  We begin with the Courant-Fischer theorem for tensors. For notational simplicity, we define $\mathbb{I}_M\define\prod\limits_{i=1}^M I_i$.

\begin{theorem}[Courant-Fischer Theorem for Tensors]\label{thm:Courant-Fischer Theorem}
Let $\mathcal{A} \in \mathbb{C}^{I_1 \times \cdots I_M \times I_1 \times \cdots I_M}$ be a Hermitian tensor with eigenvalues $\lambda_1 \geq \lambda_2 \geq \cdots \geq \lambda_{\mbox{\tiny $\mathbb{I}_M$}}$ and corresponding orthnormal eigen-tensors $\mathcal{U}_{1}, \mathcal{U}_{2}, \cdots, \mathcal{U}_{\mbox{\tiny $\mathbb{I}_M$}}$. Spaces $S_k$ and $T_k$ are obtained by spanning via eigen-tensors $\mathcal{U}_{j}$ for $j\in \left[\mathbb{I}_M \right]$ with dimensions $k$ and $\left(\mathbb{I}_M - k +1\right)$, respectively. Then, we have 
\begin{eqnarray}\label{eq:thm:Courant-Fischer Theorem}
\lambda_k &=_1& \max\limits_{\substack{S_k \subseteq \mathbb{C}^{I_1 \times \cdots I_M} \\ \dim(S_k) = k} } \min\limits_{\mathcal{X} \in S_k}\frac{\mathcal{X}^{\mathrm{T}}\star_M\mathcal{A}\star_M\mathcal{X}}{\mathcal{X}^{\mathrm{T}}\star_M\mathcal{X}}\nonumber \\
&=_2& \min\limits_{\substack{ T_k \subseteq \mathbb{C}^{I_1 \times \cdots I_M} \\ \dim(T_k) = \left(\left(\mathbb{I}_M \right) - k +1\right)} } \max\limits_{\mathbf{X} \in T_k}\frac{\mathcal{X}^{\mathrm{T}}\star_M\mathcal{A}\star_M\mathcal{X}}{\mathcal{X}^{\mathrm{T}}\star_M\mathcal{X}}.
\end{eqnarray}
\end{theorem}
\textbf{Proof:}
We will just prove the identity $=_1$ in Eq.~\eqref{eq:thm:Courant-Fischer Theorem}. The other identity $=_2$ in Eq.~\eqref{eq:thm:Courant-Fischer Theorem} can be proved similarly. 

We first show that $\lambda_k$ is achievable. As $S_k$ is the space spanned by eigen-tensors $\mathcal{U}_1,\mathcal{U}_2,\cdots,\mathcal{U}_k$. For every $\mathbf{X} \in S_k$, we can express $\mathcal{X}$ as 
\begin{eqnarray}
\mathbf{X} &=& \sum\limits_{j=1}^{k}c_j\mathcal{U}_j.
\end{eqnarray}
Then, we have
\begin{eqnarray}
\frac{\mathcal{X}^{\mathrm{T}}\star_M\mathcal{A}\star_M\mathcal{X}}{\mathcal{X}^{\mathrm{T}}\star_M\mathcal{X}}&=&\frac{\sum\limits_{j=1}^k \lambda_j c_j^2}{c_j^2} \nonumber \\
&\geq& \frac{\sum\limits_{j=1}^k \lambda_k c_j^2}{c_j^2} = \lambda_k.
\end{eqnarray}

To verify that this is the maximum eigenvalue, as $T_k$ is the space spanned by $\mathcal{U}_k,\mathcal{U}_{k+1},\cdots,\mathcal{U}_{\mbox{\tiny $\mathbb{I}_M$}}$, for any $S_k$ with dimension $k$ the intersection of $S_k$ with $T_k$ is non-empty. Then, we also have 
\begin{eqnarray}
\min\limits_{\mathcal{X} \in S_k}\frac{\mathcal{X}^{\mathrm{T}}\star_M\mathcal{A}\star_M\mathcal{X}}{\mathcal{X}^{\mathrm{T}}\star_M\mathcal{X}} &\leq& \min\limits_{\mathcal{X} \in S_k \bigcap T_k}\frac{\mathcal{X}^{\mathrm{T}}\star_M\mathcal{A}\star_M\mathcal{X}}{\mathcal{X}^{\mathrm{T}}\star_M\mathcal{X}}.
\end{eqnarray}
Any such $\mathcal{X}$ can be expressed as 
\begin{eqnarray}
\mathcal{X} &=& \sum\limits_{j=k}^{\mathbb{I}_M}c_j \mathcal{U}_j,
\end{eqnarray}
then, we have 
\begin{eqnarray}
\frac{\mathcal{X}^{\mathrm{T}}\star_M\mathcal{A}\star_M\mathcal{X}}{\mathcal{X}^{\mathrm{T}}\star_M\mathcal{X}} &=& \frac{\sum\limits_{j=k}^{\mathbb{I}_M}\lambda_j c_j^2}{c_j^2} \nonumber \\
&\leq& \frac{\sum\limits_{j=k}^{\mathbb{I}_M} \lambda_k c_j^2}{c_j^2} = \lambda_k.
\end{eqnarray}
Therefore, all subspace of $S_k$ with dimension $k$, we have 
\begin{eqnarray}
\min\limits_{\mathcal{X} \in S_k}\frac{\mathcal{X}^{\mathrm{T}}\star_M\mathcal{A}\star_M\mathcal{X}}{\mathcal{X}^{\mathrm{T}}\star_M\mathcal{X}} &\leq& \lambda_k.
\end{eqnarray}
This theorem is proved since $\lambda_k$ is achievable and is the maximum eigenvalue. 
$\hfill \Box$

We use $\mathcal{U}_{(k)} \in \mathbb{C}^{I_1 \times \cdots I_M \times I_1 \times \cdots I_M}$ to represent a projection tensor to the space with dimension $k$, i.e., $\mathcal{U}_{(k)}\star_M \mathcal{X}$ will be a tensor spanned by $k$ eigen-tensors of $\mathcal{X}$ and $\mathcal{U}_{(k)}^H \star_M \mathcal{U}_{(k)} \preceq \mathcal{I} \in \mathbb{C}^{I_1 \times \cdots I_M \times I_1 \times \cdots I_M}$. 

According to Theorem~\ref{thm:Courant-Fischer Theorem}, we can have the following lemma about the $k$-th largest eigenvalue of a random Hermitian tensor. 
\begin{lemma}\label{lma:single RT bound}
Let $\mathcal{X} \in \mathbb{C}^{I_1 \times \cdots I_M \times I_1 \times \cdots I_M}$ be a random Hermitian tensor, and $k \in [\mathbb{I}_M]$. For any $\theta \in \mathbb{R}$, we have
\begin{eqnarray}\label{eq1:lma:single RT bound}
\mathrm{Pr}\left(\lambda_k\left(\mathcal{X}\right)\geq \theta \right)&\leq&\inf\limits_{t > 0}\min\limits_{\mathcal{U}_{\left(\mathbb{I}_M - k+1\right)}}\left[e^{-t \theta} \mathbb{E} \mathrm{Tr}\exp\left(t\mathcal{U}_{\left(\mathbb{I}_M - k+1\right)}^{H}\star_M \mathcal{X} \star_M \mathcal{U}_{\left(\mathbb{I}_M - k+1\right)}\right)\right].
\end{eqnarray}
\end{lemma}
\textbf{Proof:}
Because we have
\begin{eqnarray}
\mathrm{Pr}\left(\lambda_k\left(\mathcal{X}\right)\geq \theta \right)&=_1&\mathrm{Pr}\left(e^{\lambda_k\left(t\mathcal{X}\right)}\geq e^{t\theta}\right) \nonumber \\
&\leq_2& e^{-t \theta} \mathbb{E}e^{\lambda_k\left(t\mathcal{X}\right)}  \nonumber \\
&=_3& e^{-t \theta}\mathbb{E}\exp\left(\min\limits_{\mathcal{U}_{\left(\mathbb{I}_M - k+1\right)}}\lambda_{\max}\left(t\mathcal{U}_{\left(\mathbb{I}_M - k+1\right)}^{H}\star_M \mathcal{X} \star_M \mathcal{U}_{\left(\mathbb{I}_M - k+1\right)}\right)\right) \nonumber \\
&\leq_4&  \min\limits_{\mathcal{U}_{\left(\mathbb{I}_M - k+1\right)}} \mathbb{E}\lambda_{\max}\left(\exp\left(t\mathcal{U}_{\left(\mathbb{I}_M - k+1\right)}^{H}\star_M \mathcal{X} \star_M \mathcal{U}_{\left(\mathbb{I}_M - k+1\right)}\right)\right) \nonumber \\
&\leq_5& \min\limits_{\mathcal{U}_{\left(\mathbb{I}_M - k+1\right)}}\mathbb{E}\mathrm{Tr}\left(\exp\left(t\mathcal{U}_{\left(\mathbb{I}_M - k+1\right)}^{H}\star_M \mathcal{X} \star_M \mathcal{U}_{\left(\mathbb{I}_M - k+1\right)}\right)\right),
\end{eqnarray}
where $=_1$ comes from monotonicity of the exponential function with respect to a real variable; $\leq_2$ comes from Markov inequality; $=_3$ comes from Theorem~\ref{thm:Courant-Fischer Theorem}; $\leq_4$ comes from Jensen's inequality, and $\leq_5$ comes from that the largest eigenvalue is smaller than its trace due to that the exponential of a Hemitian tensor is positive definite. This lemma is proved.
$\hfill \Box$

Below, we will present Lieb-Seiringer concavity theorem for tensors, which will be used to extend the subadditivity lemma of tensor cumulant-generating functions for the random tensor operated by projection tensors, see Lemma 4 in~\cite{chang2022convenient}.

\begin{theorem}[Lieb-Seiringer Joint Concavity for Tensors]\label{thm:Lieb-Seiringer Joint Concavity for Tensors}
Let $\mathcal{H}$ be a Hermitian tensor operating on a finite dimensional Hilbert space $\mathfrak{H}$, and $\mathcal{A}_n$ be positive-definite tensors, where $1 \leq n \leq N$. If we have 
\begin{eqnarray}\label{eq1:thm:Lieb-Seiringer Joint Concavity for Tensors}
\sum\limits_{i=1}^N \mathcal{U}^H_i \star_M \mathcal{U}_i \preceq \mathcal{I}_{\mathcal{H}},
\end{eqnarray}
where $\mathcal{U}_i \in \mathbb{C}^{I_1 \times \cdots \times I_M\times
I_1 \times \cdots \times I_M}$, and $\mathcal{I}_{\mathcal{H}}$ is the identity tensor assoicated to the Hilbert space $\mathfrak{H}$. The following mapping 
\begin{eqnarray}\label{eq2:thm:Lieb-Seiringer Joint Concavity for Tensors}
\left(\mathcal{A}_1,\cdots,\mathcal{A}_N\right) \rightarrow 
\mathrm{Tr}_{\mathfrak{H}}\exp\left(\mathcal{H}+\sum\limits_{i=1}^N \mathcal{U}^{H}_i \star_M \left(\log  \mathcal{A}_i\right)\star_M \mathcal{U}_i \right)
\end{eqnarray}
is jointly concave.
\end{theorem}
\textbf{Proof:}
We may assume that Eq.~\eqref{eq1:thm:Lieb-Seiringer Joint Concavity for Tensors} is the identity. Since if $\sum\limits_{i=1}^N \mathcal{U}^H_i \star_M \mathcal{U}_i \prec \mathcal{I}_{\mathcal{H}}$, we can add one more term $\mathcal{U}_{N+1}=\left(\mathcal{I}_{\mathfrak{H}}-\sum\limits_{i=1}^N \mathcal{U}^H_i \star_M \mathcal{U}_i \right)^{1/2}$ with $\mathcal{A}_{N+1} = \mathcal{I}_{\mathcal{H}}$.

We consider the following three transformations with respect to tensors $\mathcal{H}, \mathcal{U}_i$ and $\mathcal{A}_i$:
\begin{eqnarray}\label{eq3:thm:Lieb-Seiringer Joint Concavity for Tensors}
\tilde{\mathcal{H}}_{i,j}&=&\frac{\delta_{i,j}}{N} \mathcal{U}^H_i \star_M \mathcal{H} \star_M  \mathcal{U}_i, \nonumber \\
\tilde{\mathcal{P}}_{i,j}&=&\mathcal{U}_i \star_M \mathcal{U}^H_j , \nonumber \\
\tilde{\mathcal{A}}_{i,j}&=&\delta_{i,j}\mathcal{A}_i,
\end{eqnarray}
where $\delta_{i,j}$ is a Kronecker delta function with respect to indices $i \in [N]$ and $j \in [N]$. Note that transformed tensors $\tilde{\mathcal{H}}_{i,j}, \tilde{\mathcal{P}}_{i,j}$ and $\tilde{\mathcal{A}}_{i,j}$ are tensors operated to another Hilbert space $\mathfrak{H} \otimes \mathbb{C}^{N \times N} \define \tilde{\mathfrak{H}}$. Therefore, tensors $\tilde{\mathcal{H}}, \tilde{\mathcal{P}}$ and $\tilde{\mathcal{A}}$ are with dimensions $ \mathbb{C}^{M \times I_1 \times \cdots \times I_M\times
M \times I_1 \times \cdots \times I_M}$. Besides, we also have $\tilde{\mathcal{P}}_{i,j} = \tilde{\mathcal{P}}^H_{i,j}$ and $\tilde{\mathcal{P}}_{i,j} = \tilde{\mathcal{P}}^2_{i,j}$ since $\sum\limits_{i=1}^N \mathcal{U}^H_i \star_M \mathcal{U}_i = \mathcal{I}_{\mathcal{H}}$.

From Theorem 7 in~\cite{chang2022convenient} and any $r \in \mathbb{R}$, we have the following jointly concave map:
\begin{eqnarray}\label{eq4:thm:Lieb-Seiringer Joint Concavity for Tensors}
\left(\mathcal{A}_1,\cdots,\mathcal{A}_N\right) \rightarrow 
\mathrm{Tr}_{\tilde{\mathfrak{H}}}\exp\left(-r\left(\mathcal{I}_{\tilde{\mathfrak{H}}} - \tilde{\mathcal{P}}\right) +\tilde{\mathcal{H}} + \log \tilde{\mathcal{A}}\right),
\end{eqnarray}
where $\mathcal{I}_{\tilde{\mathfrak{H}}}$ is the identity tensor in the space $\tilde{\mathfrak{H}}$. By taking $r \rightarrow \infty$, Eq.~\eqref{eq4:thm:Lieb-Seiringer Joint Concavity for Tensors} can be reduced as
\begin{eqnarray}\label{eq5:thm:Lieb-Seiringer Joint Concavity for Tensors}
\left(\mathcal{A}_1,\cdots,\mathcal{A}_N\right) \rightarrow 
\mathrm{Tr}_{\tilde{\mathcal{P}}\tilde{\mathfrak{H}}}\exp\left(\tilde{\mathcal{P}}\star_{M+1}\left(\tilde{\mathcal{H}}+\log\tilde{\mathcal{A}}\right)\star_{M+1}\tilde{\mathcal{P}}\right).
\end{eqnarray}

By considering the map $\mathcal{V}:\mathfrak{H} \rightarrow \tilde{\mathcal{P}}\tilde{\mathfrak{H}}$ defined by
\begin{eqnarray}\label{eq6:thm:Lieb-Seiringer Joint Concavity for Tensors}
\left(\mathcal{V}\left(\mathcal{X}\right)\right)_i &=& \mathcal{U}_i \star_M \mathcal{X},
\end{eqnarray}
and the relation of $\sum\limits_{i=1}^N \mathcal{U}^H_i \star_M \mathcal{U}_i = \mathcal{I}_{\mathcal{H}}$, we know that the map $\mathcal{V}$ is a unitary map. Then, we have 
\begin{eqnarray}\label{eq7:thm:Lieb-Seiringer Joint Concavity for Tensors}
\mathcal{V}^H \star_M \left(\tilde{\mathcal{P}}\star_{M+1}\left(\tilde{\mathcal{H}}+\log\tilde{\mathcal{A}}\right)\star_{M+1}\tilde{\mathcal{P}}\right)\star_M \mathcal{V}  &=& \mathcal{H}+\sum\limits_{i=1}^N \mathcal{U}^{H}_i \star_M \left(\log  \mathcal{A}_i\right)\star_M \mathcal{U}_i. 
\end{eqnarray}
Therefore, we have 
\begin{eqnarray}\label{eq7:thm:Lieb-Seiringer Joint Concavity for Tensors}
\mathrm{Tr}_{\tilde{\mathcal{P}}\tilde{\mathcal{H}}}\exp\left(\tilde{\mathcal{P}}\star_{M+1}\left(\tilde{\mathcal{H}}+\log\tilde{\mathcal{A}}\right)\star_{M+1}\tilde{\mathcal{P}}\right)
&=& \mathrm{Tr}_{\mathcal{H}}\exp\left(\mathcal{H}+\sum\limits_{i=1}^N \mathcal{U}^{H}_i \star_M \left(\log  \mathcal{A}_i\right)\star_M \mathcal{U}_i\right)
\end{eqnarray}
This complets the proof from Eq.~\eqref{eq5:thm:Lieb-Seiringer Joint Concavity for Tensors}.
$\hfill \Box$

\textbf{Remark:} We have to provide the following comments about Theorem~\ref{thm:Lieb-Seiringer Joint Concavity for Tensors}:
\begin{enumerate}
\item Our Theorem~\ref{thm:Lieb-Seiringer Joint Concavity for Tensors} extends Theorem 3 in~\cite{lieb2005stronger} to tensors by modifying the transforms provided by Eq.~\eqref{eq3:thm:Lieb-Seiringer Joint Concavity for Tensors}, which works for both invertible and noninvertible tensors $\mathcal{U}_n$. The original transformations used by the proof in Theorem 3 require tesors $\mathcal{U}_n$ to be invertible in~\cite{lieb2005stronger}. 
\item  Proposition 3.1 utilizes Theorem 3~\cite{gittens2011tail} directly, however, the dimension of the matrix $\bm{V}$ used by Proposition 3.1 does not satisfy the requirement by Theorem 3 in~\cite{lieb2005stronger} since $\bm{V}^{H}\bm{V} = \bm{I}_k$ (identity $k \times k$ matrix), which is not the identity of the original Hilbert space $\bm{I}_n$ (The notations used here adopts from Theorem 3~\cite{gittens2011tail}). 
\item Similar to~\cite{gittens2011tail}, we also utilize Courant-Fischer theorem and Lieb-Seiringer joint concavity theorem to establish the tail bounds. However, our work extends the settings from matrices to tensors. 
\end{enumerate}

The following lemma is the subadditivity tensor cumulant-generating functions for the random tensor operated by projection tensors.
\begin{lemma}\label{lma:expectation trace exp bound}
Given $N$ indepedent Hermitian random tensors $\mathcal{X}_i \in \mathbb{C}^{I_1 \times \cdots I_M \times I_1 \times \cdots I_M}$, and $N$ deterministic Hermitian random tensors $\mathcal{A}_i \in \mathbb{C}^{I_1 \times \cdots I_M \times I_1 \times \cdots I_M}$ that satisfy
\begin{eqnarray}\label{eq1:lma:expectation trace exp bound}
\mathbb{E}\exp(\mathcal{X}_i) \preceq \exp(\mathcal{A}_i).
\end{eqnarray}
Pick $k \in [\mathbb{I}_M]$, and let $\mathcal{U}_{(k)} \in \mathbb{C}^{I_1 \times \cdots I_M \times I_1 \times \cdots I_M}$ be a projection tensor. Then, we have 
\begin{eqnarray}\label{eq2:lma:expectation trace exp bound}
\mathbb{E}\mathrm{Tr}\exp\left(\sum\limits_{i=1}^N \mathcal{U}^H_{(k)}\star_M \mathcal{X}_i \star_M \mathcal{U}_{(k)}\right)&\leq& \mathrm{Tr}\exp\left(\sum\limits_{i=1}^N \mathcal{U}^H_{(k)}\star_M \mathcal{A}_i \star_M \mathcal{U}_{(k)}\right).
\end{eqnarray}
\end{lemma}
\textbf{Proof:}
Let $\mathbb{E}_k$ be the expectation on the first $k$ random tensor from $\mathcal{X}_1$ to $\mathcal{X}_k$, then we have
\begin{eqnarray}\label{eq3:lma:expectation trace exp bound}
\lefteqn{\mathbb{E}\mathrm{Tr}\exp\left(\sum\limits_{i=1}^N \mathcal{U}^H_{(k)}\star_M \mathcal{X}_i \star_M \mathcal{U}_{(k)}\right)}\nonumber \\
&=& \mathbb{E}_0 \cdots \mathbb{E}_{N-1} \mathrm{Tr} \exp\left(\sum\limits_{i=1}^{N-1}\mathcal{U}^H_{(k)}\star_M \mathcal{X}_i \star_M \mathcal{U}_{(k)}+ \mathcal{U}^H_{(k)}\star_M \log e^{\mathcal{X}_N}\star_M \mathcal{U}_{(k)}\right) \nonumber \\
&\leq_1&  \mathbb{E}_0 \cdots \mathbb{E}_{N-2} \mathrm{Tr} \exp\left(\sum\limits_{i=1}^{N-1}\mathcal{U}^H_{(k)}\star_M \mathcal{X}_i \star_M \mathcal{U}_{(k)}+ \mathcal{U}^H_{(k)}\star_M \log \mathbb{E} e^{\mathcal{X}_N}\star_M \mathcal{U}_{(k)}\right) \nonumber \\
&\leq_2&  \mathbb{E}_0 \cdots \mathbb{E}_{N-2} \mathrm{Tr} \exp\left(\sum\limits_{i=1}^{N-1}\mathcal{U}^H_{(k)}\star_M \mathcal{X}_i \star_M \mathcal{U}_{(k)}+ \mathcal{U}^H_{(k)}\star_M \mathcal{A}_N \star_M \mathcal{U}_{(k)}\right),
\end{eqnarray}
where $=_1$ comes from Theorem~\ref{thm:Lieb-Seiringer Joint Concavity for Tensors} and Jensen's inequality, and $=_2$ comes from the condition provided by Eq.~\eqref{eq1:lma:expectation trace exp bound}. This lemma is proved by repeating the inequality relation with respect to the argument $i$ given by Eq.~\eqref{eq3:lma:expectation trace exp bound}.
$\hfill \Box$

We are ready to present the main result of this section. 
\begin{theorem}\label{thm:minmax Laplace TF}
Given $N$ indepedent Hermitian random tensors $\mathcal{X}_i \in \mathbb{C}^{I_1 \times \cdots I_M \times I_1 \times \cdots I_M}$, and $N$ deterministic Hermitian random tensors $\mathcal{A}_i \in \mathbb{C}^{I_1 \times \cdots I_M \times I_1 \times \cdots I_M}$ that satisfy
\begin{eqnarray}\label{eq1:thm:minmax Laplace TF}
\mathbb{E}\exp(t \mathcal{X}_i) \preceq \exp(f(t)\mathcal{A}_i),
\end{eqnarray}
where $f:(0, \infty) \rightarrow [0, \infty)$ and $i \in [N]$. Then, for all $\theta \in \mathbb{R}$, we have
\begin{eqnarray}\label{eq2:lma:expectation trace exp bound}
\lefteqn{\mathrm{Pr}\left(\lambda_k\left(\sum\limits_{i=1}^N\right)\geq\theta\right)
}\nonumber \\
&\leq& \inf\limits_{t >0}\min\limits_{\mathcal{U}_{\left(\mathbb{I}_M - k +1\right)}}\left[e^{-t \theta}\mathrm{Tr}\exp\left(f(t)\sum\limits_{i=1}^N \mathcal{U}_{\left(\mathbb{I}_M - k +1\right)}^H \star_M \mathcal{A}_i \star_M \mathcal{U}_{\left(\mathbb{I}_M - k +1\right)}\right)\right]
\end{eqnarray}
\end{theorem}
\textbf{Proof:}
By combining Lemma~\ref{lma:single RT bound} and Lemma~\ref{lma:expectation trace exp bound}. 
$\hfill \Box$

\section{Algebraic Connectivity of Ensemble Random Hypergraphs}\label{sec: Algebraic Connectivity of Ensemble Hypergraphs}

According to Eq.~\eqref{eq:alg def}, the value of the algebraic connectivity for a hypergraph $\mathscr{G}$ indicates the degree of connectivity within the entire hypergraph. In this section, we will consider different tail bounds for the algebraic connectivity of ensemble hypergraphs with respect to different random hypergraphs assumptions. These tail bounds are Chernoff, Bennett, and Bernstein bounds, which will be discussed in Sections~\ref{subsec: Chernoff Bounds for Algebraic Connectivity of Ensemble Hypergraphs},~\ref{subsec: Bennett Bounds for Algebraic Connectivity of Ensemble Hypergraphs}, and~\ref{subsec: Bernstein Bounds for Algebraic Connectivity of Ensemble Hypergraphs}, respectively.

\subsection{Chernoff Bounds for Algebraic Connectivity of Ensemble Hypergraphs}\label{subsec: Chernoff Bounds for Algebraic Connectivity of Ensemble Hypergraphs}

In this section, we will consider the algebraic connectivity for ensembles of random hypergraphs whose Laplacian tensors are positive semidefinite tensors.
\begin{theorem}\label{thm:Chernoff Bounds}
Given $N$ independent random hypergraphs, $\{\mathscr{G}_i\}$, whose Laplacian tensors, $\{\mathcal{L}_i \in \mathbb{R}^{m^{2M}}\}$ for $i \in [N]$, are positive semidefinite tensors. We define a positive number $\nu_{\overline{\mathscr{G}}}$ by
\begin{eqnarray}\label{eq0:thm:Chernoff Bounds}
\nu_{\overline{\mathscr{G}}} &\define& \lambda_{\frac{m!}{(m-M)!}-1} \left(\sum\limits_{i=1}^N \mathbb{E}\mathcal{L}_i\right).
\end{eqnarray}
Then, we have the following bounds for the algebraic connectivity of the ensemble hypergraph $\overline{\mathscr{G}}$: 
\begin{eqnarray}\label{eq1:thm:Chernoff Bounds}
\mathrm{Pr}\left(\alpha(\overline{\mathscr{G}}) \geq (1+\theta)\nu_{\overline{\mathscr{G}}}\right) \leq \left(m^M - \frac{m!}{(m-M)!} +2\right) \left[\frac{e^{\theta}}{(1+\theta)^{1+\theta}}\right]^{\nu_{\overline{\mathscr{G}}}},
\end{eqnarray}
where $\theta >1$; and 
\begin{eqnarray}\label{eq2:thm:Chernoff Bounds}
\mathrm{Pr}\left(\alpha(\overline{\mathscr{G}}) \leq (1+\theta)\nu_{\overline{\mathscr{G}}}\right) \leq \left(\frac{m!}{(m-M)!}-1\right) \left[\frac{e^{-\theta}}{(1-\theta)^{1-\theta}}\right]^{\nu_{\overline{\mathscr{G}}}},
\end{eqnarray}
where $0 \leq \theta <1$.
\end{theorem} 
\textbf{Proof:}
We define projection tensors $\mathcal{U}_{+}$ and $\mathcal{U}_{-}$ as 
\begin{eqnarray}
\mathcal{U}_{+}&\define& \mathcal{U}_{\left(m^M - \frac{m!}{(m-M)!} +2\right)}, \nonumber \\
\mathcal{U}_{-} &\define& \mathcal{U}_{\left(\frac{m!}{(m-M)!}-1\right)}. 
\end{eqnarray}
Then, we have 
\begin{eqnarray}\label{eq3:thm:Chernoff Bounds}
\mathrm{Pr}\left(\alpha(\overline{\mathscr{G}}) \geq (1+\theta)\nu_{\overline{\mathscr{G}}}\right) 
&=_1& \mathrm{Pr}\left(\lambda_{\frac{m!}{(m-M)!}-1}\left(\sum\limits_{i=1}^N \mathcal{L}_i\right)\geq (1+\theta)\nu_{\overline{\mathscr{G}}}\right) \nonumber\\
&\leq_2& 
\inf\limits_{t > 0} e^{-t(1+\theta)\nu_{\overline{\mathscr{G}}}} \mathrm{Tr}\exp\left((e^t - 1)\sum\limits_{i=1}^N \mathcal{U}_{+}^H \star_M \mathbb{E}\mathcal{L}_i\star_M \mathcal{U}_{+}\right)\nonumber\\
&\leq_3&\inf\limits_{t > 0} e^{-t(1+\theta)\nu_{\overline{\mathscr{G}}}}\left(m^M - \frac{m!}{(m-M)!} +2\right)\nonumber \\
&& \cdot \lambda_{\max}\left(\exp\left((e^t - 1)\sum\limits_{i=1}^N \mathcal{U}_{+}^H \star_M \mathbb{E}\mathcal{L}_i\star_M \mathcal{U}_{+}\right)\right) \nonumber \\
&=_4& \left(m^M - \frac{m!}{(m-M)!} +2\right)\inf\limits_{t > 0}\exp\left[\left(e^t -1 - t (1 +\theta)\right) \nu_{\overline{\mathscr{G}}}\right],
\end{eqnarray}
where $=_1$ comes from the fact that the algebraic connectivity $\alpha(\overline{\mathscr{G}})$ is the $\left(\frac{m!}{(m-M)!}-1\right)$-th largest eigenvalue of $\sum\limits_{i=1}^N \mathcal{L}_i$ from Lemma 3.1 in~\cite{gu2023even}, $\leq_2$ comes from Theorem~\ref{thm:minmax Laplace TF} and Lemma 6 in~\cite{chang2022convenient}, $\leq_3$ comes from the fact that the trace can be bounded by the maximum eigenvalue and the dimensions size, and $=_4$ comes from Theorem~\ref{thm:Courant-Fischer Theorem} and Eq.~\eqref{eq0:thm:Chernoff Bounds}.  Eq.~\eqref{eq1:thm:Chernoff Bounds} is proved by applying the minimizer $t=\log(1+\theta)$ at the term $\exp\left[\left(e^t -1 - t (1 +\theta)\right) \nu_{\overline{\mathscr{G}}}\right]$. 

Let us prove Eq.~\eqref{eq2:thm:Chernoff Bounds}. We have 
\begin{eqnarray}\label{eq4:thm:Chernoff Bounds}
\mathrm{Pr}\left(\alpha(\overline{\mathscr{G}}) \leq (1-\theta)\nu_{\overline{\mathscr{G}}}\right) 
&=_1& \mathrm{Pr}\left(\lambda_{\frac{m!}{(m-M)!}-1}\left(\sum\limits_{i=1}^N \mathcal{L}_i\right)\leq (1-\theta)\nu_{\overline{\mathscr{G}}}\right) \nonumber\\
&=_2& \mathrm{Pr}\left(\lambda_{m^M - \frac{m!}{(m-M)!}+2}\left(-\sum\limits_{i=1}^N \mathcal{L}_i\right)\geq -(1-\theta)\nu_{\overline{\mathscr{G}}}\right) \nonumber\\
&\leq_3& 
\inf\limits_{t > 0} e^{t(1-\theta)\nu_{\overline{\mathscr{G}}}} \mathrm{Tr}\exp\left((1-e^{-t })\sum\limits_{i=1}^N \mathcal{U}_{-}^H \star_M (-\mathbb{E}\mathcal{L}_i)\star_M \mathcal{U}_{-}\right)\nonumber\\
&\leq_4&\inf\limits_{t > 0} e^{t(1-\theta)\nu_{\overline{\mathscr{G}}}} \left(\frac{m!}{(m-M)!}-1\right)\nonumber \\
&& \cdot \lambda_{\min}\left(\exp\left((e^{-t} - 1)\sum\limits_{i=1}^N \mathcal{U}_{-}^H \star_M \mathbb{E}\mathcal{L}_i\star_M \mathcal{U}_{-}\right)\right) \nonumber \\
&=_5& \left(\frac{m!}{(m-M)!}-1\right)\inf\limits_{t > 0}\exp\left[\left(t(1-\theta)-1+e^{-t}\right) \nu_{\overline{\mathscr{G}}}\right],
\end{eqnarray}
where $=_1$ comes from the fact that the algebraic connectivity $\alpha(\overline{\mathscr{G}})$ is the $\left(\frac{m!}{(m-M)!}-1\right)$-th largest eigenvalue of $\sum\limits_{i=1}^N \mathcal{L}_i$ from Lemma 3.1 in~\cite{gu2023even}, $=_2$ comes from the eigenvalues duality fact that \\
$\lambda_{\left(\frac{m!}{(m-M)!}-1\right)}\left(-\sum\limits_{i=1}^N \mathcal{L}_i\right) = - \lambda_{\left(m^M - \frac{m!}{(m-M)!}+2\right)}\left(\sum\limits_{i=1}^N \mathcal{L}_i\right)$, $\leq_3$ comes from Theorem~\ref{thm:minmax Laplace TF} and Lemma 6 in~\cite{chang2022convenient}, $\leq_4$ comes from the fact that the trace can be bounded by the maximum eigenvalue and the dimensions size, and $=_5$ comes from Theorem~\ref{thm:Courant-Fischer Theorem} and Eq.~\eqref{eq0:thm:Chernoff Bounds}. Then, Eq.~\eqref{eq2:thm:Chernoff Bounds} is proved by applying the minimizer $t=-\log(1-\theta)$ at the term $\exp\left[\left(t(1-\theta)-1+e^{-t}\right) \nu_{\overline{\mathscr{G}}}\right]$. 
$\hfill \Box$

\subsection{Bennett Bounds for Algebraic Connectivity of Ensemble Hypergraphs}\label{subsec: Bennett Bounds for Algebraic Connectivity of Ensemble Hypergraphs}

In this section, we will consider the algebraic connectivity for ensemble of random hypergraphs whose  Laplacian tensors are Hermitian.

\begin{theorem}\label{thm:Bennett Bounds}
Given $N$ independent random hypergraphs, $\{\mathscr{G}_i\}$, whose Laplacian tensors, $\{\mathcal{L}_i \in \mathbb{R}^{m^{2M}}\}$ for $i \in [N]$, are Hermitian tensors with $\mathbb{E}\mathcal{L}_i = \mathcal{O}$ and $\lambda_{\max}(\mathcal{L}_i) \leq 1$ almost surely. We define a positive number $\sigma^2_{\overline{\mathscr{G}}}$ by
\begin{eqnarray}\label{eq0:thm:Bennett Bounds}
\sigma^2_{\overline{\mathscr{G}}} &\define& \lambda_{\frac{m!}{(m-M)!}-1} \left(\sum\limits_{i=1}^N \mathbb{E}\mathcal{L}^2_i\right).
\end{eqnarray}
Then, we have the following bounds for the algebraic connectivity of the ensemble hypergraph $\overline{\mathscr{G}}$: 
\begin{eqnarray}\label{eq1:thm:Bennett Bounds}
\mathrm{Pr}\left(\alpha(\overline{\mathscr{G}}) \geq \theta \right)&\leq&\left(m^M - \frac{m!}{(m-M)!} +2\right)\exp(\theta)\left(1+\frac{\theta}{\sigma^2_{\overline{\mathscr{G}}}}\right)^{-(\theta+\sigma^2_{\overline{\mathscr{G}}})}.
\end{eqnarray}
\end{theorem}
\textbf{Proof:}
We define projection tensors $\mathcal{U}_{+}$as 
\begin{eqnarray}
\mathcal{U}_{+}&\define& \mathcal{U}_{\left(m^M - \frac{m!}{(m-M)!} +2\right)}.
\end{eqnarray}
Then, we have
\begin{eqnarray}\label{eq2:thm:Bennett Bounds}
\mathrm{Pr}\left(\alpha(\overline{\mathscr{G}}) \geq \theta \right)&=_1& \mathrm{Pr}\left(\lambda_{\frac{m!}{(m-M)!}-1}\left(\sum\limits_{i=1}^N \mathcal{L}_i\right)\geq \theta\right) \nonumber\\
&\leq_2&\inf\limits_{t >0}e^{-t\theta}\mathrm{Tr}\exp\left(\left(e^t - t -1\right)\sum\limits_{i=1}^N \mathcal{U}_{+}^H \star_M \mathbb{E}\mathcal{L}^2_i\star_M \mathcal{U}_{+}\right)\nonumber\\
&\leq_3&\left(m^M - \frac{m!}{(m-M)!} +2\right)\nonumber \\
&& \cdot \inf\limits_{t > 0}\left[e^{-t\theta}\lambda_{\max}\left(\exp\left((e^t -t -1)\sum\limits_{i=1}^N \mathcal{U}_{+}^H \star_M \mathbb{E}\mathcal{L}^2_i\star_M \mathcal{U}_{+}\right)\right)\right]\nonumber \\
&=_4& \left(m^M - \frac{m!}{(m-M)!} +2\right)\inf\limits_{t>0}e^{(e^t -t -1) \sigma^2_{\overline{\mathscr{G}}}- t \theta}
\end{eqnarray}
where $=_1$ comes from the fact that the algebraic connectivity $\alpha(\overline{\mathscr{G}})$ is the $\left(\frac{m!}{(m-M)!}-1\right)$-th largest eigenvalue of $\sum\limits_{i=1}^N \mathcal{L}_i$ from Lemma 3.1 in~\cite{gu2023even}, $\leq_2$ comes from Theorem~\ref{thm:minmax Laplace TF} and Lemma 7 in~\cite{chang2022convenient}, $\leq_3$ comes from the fact that the trace can be bounded by the maximum eigenvalue and the dimensions size, and $=_4$ comes from Theorem~\ref{thm:Courant-Fischer Theorem} and Eq.~\eqref{eq0:thm:Bennett Bounds}.  Eq.~\eqref{eq1:thm:Bennett Bounds} is proved by applying the minimizer $t=\log(1+\theta/\sigma^2_{\overline{\mathscr{G}}})$ at the term $(e^t -t -1) \sigma^2_{\overline{\mathscr{G}}}- t \theta$. 
$\hfill \Box$

\subsection{Bernstein Bounds for Algebraic Connectivity of Ensemble Hypergraphs}\label{subsec: Bernstein Bounds for Algebraic Connectivity of Ensemble Hypergraphs}

In this section, we will consider the algebraic connectivity for ensembles of random hypergraphs whose  Laplacian tensors are Hermitian with subexponential growth rate. 

\begin{theorem}\label{thm:Bernstein Bounds}
Given $N$ independent random hypergraphs, $\{\mathscr{G}_i\}$, whose Laplacian tensors, $\{\mathcal{L}_i \in \mathbb{R}^{m^{2M}}\}$ for $i \in [N]$, are Hermitian tensors. We assume that these Laplacian random tensors satisfy the following subexponential growth rate as:
\begin{eqnarray}\label{eq1:thm:Bernstein Bounds}
\mathbb{E}(\mathcal{L}_i^p) \preceq \frac{p! \mathcal{A}_i^2}{2},
\end{eqnarray}
where $p=2,3,4,\cdots,$, $\mathbb{E}(\mathcal{L}_i) = \mathcal{O}$ and $\mathcal{A}_i$ are positive-definite tensors. We also define $\sigma^2_{\overline{\mathscr{G}}}$ by 
\begin{eqnarray}\label{eq0:thm:Bernstein Bounds}
\sigma^2_{\overline{\mathscr{G}}} &\define& \lambda_{\frac{m!}{(m-M)!}-1} \left(\sum\limits_{i=1}^N \mathbb{E}\mathcal{L}^2_i\right).
\end{eqnarray}
Then, we have
\begin{eqnarray}\label{eq1:thm:Bernstein Bounds}
\mathrm{Pr}\left(\alpha(\overline{\mathscr{G}}) \geq \theta \right)
&\leq& \left(m^M - \frac{m!}{(m-M)!} +2\right)e^{-\frac{\theta^2}{2\left(\theta+\sigma^2_{\overline{\mathscr{G}}}\right)}}.
\end{eqnarray}
\end{theorem}
\textbf{Proof:}
We define projection tensors $\mathcal{U}_{+}$as 
\begin{eqnarray}
\mathcal{U}_{+}&\define& \mathcal{U}_{\left(m^M - \frac{m!}{(m-M)!} +2\right)}.
\end{eqnarray}
Then, we have
\begin{eqnarray}\label{eq2:thm:Bernstein Bounds}
\mathrm{Pr}\left(\alpha(\overline{\mathscr{G}}) \geq \theta \right)&=_1& \mathrm{Pr}\left(\lambda_{\frac{m!}{(m-M)!}-1}\left(\sum\limits_{i=1}^N \mathcal{L}_i\right)\geq \theta\right) \nonumber\\
&\leq_2&\inf\limits_{0<t<1}e^{-t\theta}\mathrm{Tr}\exp\left(\left(\frac{t^2}{2(1-t)}\right)\sum\limits_{i=1}^N \mathcal{U}_{+}^H \star_M \mathbb{E}\mathcal{L}^2_i\star_M \mathcal{U}_{+}\right)\nonumber\\
&\leq_3&\left(m^M - \frac{m!}{(m-M)!} +2\right)\nonumber \\
&& \cdot \inf\limits_{0<t<1}\left[e^{-t\theta}\lambda_{\max}\left(\exp\left(\left(\frac{t^2}{2(1-t)}\right)\sum\limits_{i=1}^N \mathcal{U}_{+}^H \star_M \mathbb{E}\mathcal{L}^2_i\star_M \mathcal{U}_{+}\right)\right)\right]\nonumber \\
&=_4& \left(m^M - \frac{m!}{(m-M)!} +2\right)\inf\limits_{0<t<1}e^{\left(\frac{t^2}{2(1-t)}\right) \sigma^2_{\overline{\mathscr{G}}}- t \theta},
\end{eqnarray}
where $=_1$ comes from the fact that the algebraic connectivity $\alpha(\overline{\mathscr{G}})$ is the $\left(\frac{m!}{(m-M)!}-1\right)$-th largest eigenvalue of $\sum\limits_{i=1}^N \mathcal{L}_i$ from Lemma 3.1 in~\cite{gu2023even}, $\leq_2$ comes from Theorem~\ref{thm:minmax Laplace TF} and Lemma 8 in~\cite{chang2022convenient}, $\leq_3$ comes from the fact that the trace can be bounded by the maximum eigenvalue and the dimensions size, and $=_4$ comes from Theorem~\ref{thm:Courant-Fischer Theorem} and Eq.~\eqref{eq0:thm:Bernstein Bounds}.  Eq.~\eqref{eq1:thm:Bernstein Bounds} is proved by applying the minimizer $t=\frac{\theta}{\theta+\sigma^2_{\overline{\mathscr{G}}}}$ at the term $\left(\frac{t^2}{2(1-t)}\right) \sigma^2_{\overline{\mathscr{G}}}- t \theta$. 
$\hfill \Box$

\bibliographystyle{IEEETran}
\bibliography{TenConHypergraph_Bib}

\end{document}